\documentclass[a4paper,review,authoryear]{elsarticle}

\usepackage{amsfonts}
\usepackage{amssymb}
\usepackage{enumerate}
\usepackage{amsmath}
\usepackage{graphicx}
\usepackage{theorem}
\usepackage{array}
\usepackage{natbib}
\usepackage{lscape}
\usepackage{rotating}
\usepackage{url}

\newtheorem{theorem}{Theorem}

\newtheorem{corollary}{Corollary}

{\theorembodyfont{\normalfont}
\newtheorem{proofprop}{Proof of Proposition}
}
{\theorembodyfont{\normalfont}
\newtheorem{prooftheorem}{Proof of Theorem}
}
{\theorembodyfont{\normalfont}

}
{\theorembodyfont{\normalfont}

}
{\theorembodyfont{\normalfont}
\newtheorem{remark}{Remark}
}
{\theorembodyfont{\normalfont}

}
{\theorembodyfont{\normalfont}
\newtheorem{prop}{Proposition}
}

\begin{document}

\title{New Approach to Derive the Value Function of a Firm with Exit Option}
\tnotetext[t1]{The research of the authors was partially supported by the SANAF project, with project reference number UTA\_CMU/MAT/0006/2009.}

\author[la]{M. Guerra}
\author[li]{C. Oliveira}
\author[li]{C. Nunes\corref{cor1}}

\cortext[cor1]{Corresponding author, Email: cnunes@math.ist.utl.pt, Phone: +351 965694036}

\address[la]{Department of Mathematics and CEMAPRE, Instituto Superior de Economia e Gest\~ao (ISEG), Rua do Quelhas n. 6, 1200-781 Lisboa, Portugal}
\address[li]{Department of Mathematics and CEMAT, Instituto Superior T{\'e}cnico, Av. Rovisco Pais, 1049-001 Lisboa, Portugal}

\begin{abstract}

In this paper we propose a new way of proving the value of a firm that is currently producing a certain product and faces the option to exit the market. 

The problem of optimal exiting  is an optimal stopping problem, that can be solved using 
the dynamic programming principle.  This approach leads to a partial differential equation, called the Hamilton-Jacobi-Bellman equation. 

This is a free-boundary problem, and therefore, we propose an approximation for the original model. 
We prove the convergence of the solution of the approximated problem to the original one and finally, using the Implicit Function Theorem, we obtain this solution.
\end{abstract}

\maketitle

\section{Introduction}

In this paper we propose a new way to prove the value function of a company that faces the option to exit from the market. 

This is a problem often adressed in the literature of real options, specially after the pionner work of \cite{DixitPindyck}. Since then, many irreversible investment problems have been widely studied in economic literature, but also in mathematical journals, due to the challenging questions that such problems raise.

In real options models the companies (that produce the goods) make decisions concerning labour levels and capital investment.  These decisions share three important characteristics. First, the  decision is partially or completely irreversible and involves some sunk costs. Second, there is uncertainty over the future rewards from the investment. Third, there is some flexibility about the timing of the decision. One can postpone the decision to get more information about the future.

There are many possible features in such an investment problem, namely the dynamics of the underlying process (finite versus infinite horizon), the possible options (exit, investment, suspension, others), the number of possible investments (only one investment versus multiple investments), type of costs (irreversible versus reversible). The list of publications is vast. For example, \cite{BobtcheffVilleneuve_2010} analyze a model of irreversible investment with two sources of uncertainty, applying their methodology to power generation under uncertainty; \cite{Wu20121241} use options game evaluation framework to study a company in the TFT-LCD industry. \cite{LokkaZervos_2011} consider an investment project where the capacity can be expanded irreversibly over time; other models related with capacity expansion can be found, for example, in \cite{AbelEberly_1996, ChiarollaHaussmann_2005, VanMieghem}. Related with the technological adoption, we refer to \cite{Farzinetal_1998,Huisman,Hagspieletal_2011b}.

In the mathematical economics literature some reversible investment problems are formulated as singular stochastic control problems. For example, \cite{GuoPham_2005} and \cite{MerhiZervos_2007} address such problem in infinite horizon, whereas \cite{HamadeneJeanblank_2007} address in finite horizon.

In the paper we study the Hamilton-Jacobi-Bellman (HJB, for short) equation (that can be obtained by the Dynamical Programming Principle and the infinitesimal generator), which turns out to be a free-boundary problem. The novelty that we propose in this work is the methodology that we use in order to circumvent the insufficient number of boundary conditions. We propose a truncation method, where natural boundary conditions exist, and such that the solution of this problem converges to the solution of the original problem.

Although the model that we propose in this paper is quite simple (dynamics of the problem is described by a geometric Brownian motion, we assume infinite horizon, and there is just one investment option - the exit of the market, and irreversible costs), this new way of deriving the optimal value function of the firm, and the optimal stopping time (the time at which the company decides, optimally, to exit the market) that we propose can also be extended to other more realistic and complex situations.

The rest of the paper is organized in the following way: in Section 2 we present the model, as well as the the conditions assumed throughout the work. In Section 3 we introduce the truncated problem, which will be an approximation of the original model, and we provide important results concerning the Hamilton-Jacobi-Bellman equation associated to this new problem. In Section 4 we prove the convergence of the solution of the truncated problem to the original one and, finally, we provide its solution, notably the optimal value function and the optimal investment decision, providing the value of the demand that triggers the exit decision. Finally  in Section 5 we present some conclusions.

\section{Model}
In this section we present the model and assumptions that we consider along the paper.

The firm currently produces an established product. The demand of this product is modelled by the unidensional stochastic process $X=\{X(t), t\geq 0\}$,  defined on a complete filtered space $((\Omega, \{{\cal F}_t\}_{t\geq 0}, {\cal P})$, with dynamics
\begin{align}
&dX(t)=\alpha X(t)dt+\sigma X(t)dW(t)\\
&X(0)=x
\end{align}
where $W=\{W(t), t \geq 0\}$ is a Brownian motion. Thus it follows that $X$ is a geometric Brownian motion, with drift $\mu$ and volatility $\sigma$. 

The price of the product at time $t$, hereby denoted by $p$, is determined by the demand and the quantity, denoted by $q$, by the following inverse demand function: 
$$p(x)=\gamma x - q,$$ 
where $x$ is the present demand level of the produce. We assume that the firm's profit, here denoted by $\Pi$, is given by 
$$\Pi(x)=p(x)\times q-K,$$ 
where $K$ is a fixed cost. 

In this paper we assume that the firm capacity is flexible, and produces optimally, in order to maximize its profit. So trivial calculations lead to the following optimal quantity, given that at the present time the demand level is $x$:
$$q(x)=\frac{\gamma x}{2}$$ 
and therefore the maximum profit is 
$$\Pi(x)=\frac{\gamma^2x^2}{4}-K$$ 

At any time, the firm may decide to exit the market, and this decision is irreversible. If the firm decides to exit the market, the firm has to pay a sunk cost, which we denote by $I\in\mathbb{R}$, independently of the actual demand level at the exiting time. when $I\geq 0$, we can interpret $I$ as a liquidation  cost, whereas otherwise it is a salvage value.

We let ${\cal S}$ denote the set of all ${\cal F}_t$-adapted Markov Times. 

The decision that the company faces is to decide optimally when it should exit the market. Thus we want to solve the following optimal stopping problem: 
\begin{align}
V(x)&=\sup_{\tau \in {\cal S}} J(x,\tau),  \label{optimal0}\\
J(x,\tau)&= E\left[\int_0^{\tau}e^{-r s}\Pi(X(s))ds-e^{-r\tau}I | X(0)=x\right], \label{optimal00}
\end{align}
where $r$ is the (constant) interest rate. In Equation (\ref{optimal00}), $\tau$ denotes the time that the company decides to exit the market, so that $J(x,\tau)$ denotes the value of the company given that the actual demand level is $x$ (initial demand level) and that it will exit the market at time $\tau$. Moreover $V(x)$ is the maximum value of the firm.

Notice that trivial calculations show that one may incorporate the sunk cost $I$ in the cost $K$.
$$
e^{-r\tau}I=I+(e^{-r\tau}-1)I=I-\int_0^{\tau}e^{-rs}\tau Ids.
$$ 
So, henceforward we only discuss the problem
\begin{align}
V(x)&=\sup_{\tau \in {\cal S}} J(x,\tau)  \label{optimal}\\
J(x,\tau)&= E\left[\int_0^{\tau}e^{-r s}\Pi(X(s))ds| X(0)=x\right]. \nonumber
\end{align}
We note that if $K\leq 0$, then we have a trivial problem, because in that case $\Pi(x)$ is always positive for all $x\in [0,+\infty[$ and therefore we do not stop, meaning that the company will  never exit the problem.   Thus during the rest of the paper, we assume that $K$ is strictly positive.

For ease the notation, in this paper we always assume the following notation: $E[\cdot]\equiv E[\cdot|X(0)=x].$

Furthermore, we also assume that $E\left[\int^{+\infty}_0e^{-r t}\left|\Pi(X(t))\right|dt\right]<+\infty$, see for example \cite{oksendal}. We note that as the profit function $\Pi(.)$ is lower-bounded,  we have that
\begin{equation*}
\label{condition}
E\left[\int^{+\infty}_0e^{-r t}\Pi^-(X(t))dt\right]<+\infty,   
\end{equation*}
whereby $E\left[\int^{+\infty}_0e^{-r t}\Pi(X(t))dt\right]$ is well defined. When $E\left[\int^{+\infty}_0e^{-r t}\Pi(X(t))dt\right]=+\infty$ the problem is trivial, leading to $\tau=\infty$. The next Proposition gives sufficient conditions for we have
$$E\left[\int^{+\infty}_0e^{-r t}\Pi(X(t))dt\right]=+\infty$$.
\begin{prop} 
\begin{itemize}
\item[(a)]If there  are $a,b\in\mathbb{R}$ such that 
$$\Pi(x)\leq a+bx^{\beta} \quad, \forall x>0\quad \text{and}\quad -r+\frac{{\sigma}^2}{2}(\beta-1)\beta+\beta\alpha<0,$$ then $E[\int_0^{+\infty}e^{-rs}\Pi(X(s))ds]<+\infty$. 
\item[(b)]If there  are $a,b\in\mathbb{R}$ such that $$\Pi(x)\geq a+bx^{\beta}, \quad \forall x>0 \quad\text{and}\quad -r+\frac{{\sigma}^2}{2}(\beta-1)\beta+\beta\alpha\geq 0,$$ then $E[\int_0^{+\infty}e^{-rs}\Pi(X(s))ds]=+\infty$.  
\end{itemize}
\end{prop}
\begin{proofprop}
Using Fubini's theorem we guarantee that \\$E[\int_0^{+\infty}e^{-rs}\Pi(X(s))ds]=\int_0^{+\infty}e^{-rs}E[\Pi(X(s))]ds$. Moreover 
\begin{align*}
E[a+b\left(X(s)\right)^{\beta}]=&E\left[a+bx_0^{\beta}\exp\Big(\beta(\alpha-\frac{{\sigma}^2}{2})s+\beta\sigma W(s)\Big)\right]\\
=&a+bx_0^{\beta}\exp\Big(\beta\alpha s+\frac{{\sigma}^2}{2}(\beta-1)\beta s\Big)
\end{align*}
Since $x_0$ is finite, we have (a) and (b). \qed
\end{proofprop}
\begin{remark}
Our function $\Pi(.)$ is polynomial on the demand level, so we have strict equality in the previous proposition with $\beta=2$. Therefore, we have two situations:
\begin{itemize}
\item[(a)] If $-r+\sigma^2+2\alpha\geq 0$, then  $E[\int_0^{+\infty}e^{-rs}\Pi(X(s))ds]=+\infty$. Therefore the problem is trivial and  $\tau=+\infty$.
\item[(b)] If $-r+\sigma^2+2\alpha< 0$, then $E[\int_0^{+\infty}e^{-rs}\Pi(X(s))ds]<+\infty $,   and therefore $\tau=+\infty$ is not {\em a priori} optimal.
\end{itemize}  
Henceforward, we assume the condition 
\begin{equation}
\label{condition1}
r>\sigma^2+2\alpha.
\end{equation}
\end{remark}

\section{The Truncated Problem}
The problem of optimal decision concerning the exit of the market is an optimal stopping problem. This kind of problems can be solved using the so-called HJB,  derived using the infinitesimal generator of the process (here assumed to be a geometric Brownian motion), and the dynamical programming principle. 

When using the HJB, we end up with a differential equation, and thus in order to solve the problem we need to find the solution of this differential equation. There are natural boundary conditions, but in this class of problems the number of boundary conditions is not enough to assure an unique solution. So we end up with a free-boundary problem \cite{Shiryaev}.

In this section we present a new method that we propose to use in order to circumvent the free-boundary problem, using a truncated problem (an approximation of the original one). This truncated problem is no-longer a free-boundary problem, and also its solution converges to the solution of the original problem. 

In the next sub-section we present in more detail the truncated problem.

\subsection{The Model}

For each $C\in ]0,+\infty[$, let $\nu_C=\inf\{t:X(t) \geq C\}$ denote the first time the demand level hits the level $C$. For such $C$, we consider the truncated process $\{X^C(t), t \geq 0\}$, with:
$$X^C(t)=
\begin{cases}
X(t) \quad\text{if}\quad t < \nu^C\\
0 \quad\text{if}\quad t \geq\nu^C
\end{cases}. 
$$
Here, $\nu_C$ is a Markov time and $X^C$ is an adapted process. Thus $X^C$ is path-wise equal to $X$ until the time it first hits level $C$; after this time, the process stays freezed in level 0. Thus for this new problem, if the firm did not exit the market until time $\nu^C$, then it will do it at this time.

For such problem we look for the the solution, hereby denoted by $V_C(.)$, such that  
\begin{align}
\label{modified_problem}
V_C(x)\equiv &\sup_{\tau \in {\cal S}} E\left[\int_0^{\tau}e^{-r s}\Pi(X^c(s))ds\right]\\
\equiv&\sup_{\tau \in {\cal S}}J_C(x,\tau)\nonumber 
\end{align}

\subsection{The Hamilton-Jacobi-Bellman Equation}

The HJB equation for the truncated problem (\ref{modified_problem}) is the following:
\begin{align}
\label{HJBtruncated0}
&\min\left\{rV_C(x)-\alpha xV_C'(x)-\frac{1}{2}\sigma^2 x^2V_C''(x)-\frac{\gamma^2x^2}{4}+K,V_C(x)\right\}=0 \quad \forall x \geq 0\\
\label{HJBtruncated}
&V_C(C)=0 \\
&V_C(0)=0 \label{boundary2}
\end{align}
where the bounday conditions (\ref{HJBtruncated}, \ref{boundary2}) follow in view of the definition of the truncated problem. Remark that this is no longer a free-boundary problem.

We start solving the following ordinary differential equation (ODE):   
\begin{equation}
\label{edo}
r v(x)- \alpha xv'(x)-\frac{1}{2}\sigma^2 x^2v''(x)- \frac{\gamma^2x^2}{4}+K=0
\end{equation}
corresponding to the value of the firm at the continuation region (i.e., where it is optimal to continue producing).

Considering the change of variable $x=e^t$, we end up with the following linear second order ODE, $$r y+ (\frac{1}{2}\sigma^2-\alpha) y'-\frac{1}{2}\sigma^2 y''= \frac{\gamma^2e^{2t}}{4}-K.$$
The solution of the homogeneous equation is $y_h=C_1e^{D_1t}+C_2e^{D_2t}$, where $D_1$ and $D_2$ are the positive and negative solution, respectively, of the equation $-\frac{1}{2}\sigma^2D^2+(\frac{1}{2}\sigma^2-\alpha) D+r=0$. Using the condition (\ref{condition1}), it is trivial to prove that $D_1>2$.

A particular solution for the non-homogeneous equation is $y_p=\frac{\gamma^2}{4(r-2\alpha-\sigma^2)}e^{2t}-\frac{K}{r}$. Therefore the solution of the equation (\ref{edo}) is given by: 
\begin{align}
\label{EDOsolution}
&v(x)=A_1x^{D_1}+A_2x^{D_2}+\frac{\gamma^2}{4(r-2\alpha-\sigma^2)}x^{2}-\frac{K}{r}.
\end{align}
Next we verify that (\ref{EDOsolution}) is indeed a solution of the problem (\ref{HJBtruncated}). For such purpose, we use the verification theorem (\ref{verification}) (see Appendix (\ref{A})). In order to use this theorem we need to assume that the solution to this problem is, at least, $C^1(\mathbb{R})$ with absolutely continuous derivative. Suppose that there is, $v_C$, a solution of (\ref{HJBtruncated}) that satisfies:
\begin{equation}
\label{fitconditions}
\begin{cases}
&v_C(x_C^*)=0 \quad \text{and} \quad v_C'(x_C^*)=0, \quad \text{(fit condition)}\\
&v_C(C)=0 \quad \quad \quad \quad \quad \quad \quad \text{(terminal condition)}.
\end{cases}
\end{equation}
where $x_C^*\in[0,+\infty[$, for each  $C\in ]0,+\infty[$.

We remark that if the solution of (\ref{edo}) is such that $v_c(0)=v_c(C)=0$, then $A_2=0, K=0$ and we end up in the trivial case. So let us assume that we are not in the trivial case.


We are going to prove that $V(x)=\max(v_C,0)$ and we guess that the continuation region for the truncated problem has the form $D_c=]x^*_c,C[$. Notice that HJB equation (\ref{HJBtruncated0}) are verified if $v_C(x)>0 ~~ \forall x\in ]x^*_C, C[$ and $r\times 0 -\alpha x \times 0-\frac{\gamma^2x^2}{4}+K\geq 0$ when $x\in ]0,x^*_C[$. Then, if the following conditions are verified:
\begin{align*}
&v_C(x)\geq 0, \quad \text{if} \quad x\geq x^*_C\\
&x^*_C\leq \frac{2K^{\frac{1}{2}}}{\gamma},
\end{align*}
the verification theorem (\ref{verification}) guarantees that $V_C(x)$, defined as
\begin{equation*}
V_C(x)=
\begin{cases}
v_C(x), \quad x\geq x^*_C\\
0, \quad x<x^*_C.
\end{cases}
\end{equation*}
is the solution of truncated problem.

\begin{prop}
The function, $v_c(x)$, defined by (\ref{EDOsolution}) with boundary conditions (\ref{fitconditions}) is positive when $x\in (x^*_c,C)$.
\end{prop}
\begin{proofprop}
Suppose that there is $\tilde{x}\in (x^*_c,C)$ such that $v_c(\tilde{x})\leq 0$. By the continuity of the function $v_c(.)$ is trivial to argue tha there exist $x_1\in  (x^*_c,C)$ such that $v_c(x_1)=\min\{v_c(x):x\in (x^*,c)\}\leq 0$. Naturally, $v_c'(x_1)=0$ and $v_c''(x_1)\geq 0$.
So,
$$
rv_c(x_1)-\frac{1}{2}\sigma^2x_1^2v_c''(x_1)=-\left(K-\frac{\gamma^2 x_1^2}{4}\right).
$$
This implies that $K-\frac{\gamma^2 x_1^2}{4}\geq 0$ and consequently $K-\frac{\gamma^2 {x^*_c}^2}{4}> 0$.
The first conclusion is, if there is a minimum point of $v_c$, first there is a maximum point of $v_c$, because $v_c(x_c^*)=v_c'(x_c^*)=0$ and $\frac{\sigma^2}{2}(x^*_c)^2v_c''(x^*_c)=K-\frac{\gamma^2(x^*_c)^2}{4}>0$, by smoothness of $v_c$ at $x^*_c$.

Therefore, let be $x_2$ such point, with $x^*_c<x_2<x_1<C$. 
Then:
\begin{align*}
&v_c(x_2)\leq 0, \quad  v_c'(x_2)=0 \quad \text{and} \quad v_c''(x_2)\leq 0.
\end{align*}
and thus $K-\frac{\gamma^2x_2^2}{4}<0$ and $K-\frac{\gamma^2x_1^2}{4}>0$ and, consequently, $x_1<x_2$, which is a contradiction resultant of our initial hypothesis. Therefore $v_c$ is positive in this region. \qed
\end{proofprop}

\begin{corollary}
The value $K-\frac{\gamma^2 {(x^*_c)}^2}{4}$ is non-negative.
\end{corollary}
The previous proposition and corollary guarantee that the conditions of the verification theorem (\ref{verification}) are verified, and so $V_c(x)$ is indeed the solution of the truncated problem.

\section{Solution of the Original Problem}
In this section we find the solution of the original problem. For that we prove that the solution of the truncated problem converges to the original one.
\begin{theorem}
The solution of the original problem, (\ref{optimal}), $V(x)$, is given by
$$
\begin{cases}
v_{\infty}(x), \quad x\geq x^*_{\infty }\\
0, \quad x<x^*_{\infty},
\end{cases},
$$
where $v_{\infty}$ is given by (\ref{EDOsolution}) with parameters $x^*_{\infty}, A^{\infty}_1, A^{\infty}_2$ given by
\begin{equation}
\label{solution_infty}
\begin{cases}
&A_1^{\infty}=0\\
&x^*_{\infty}=\left(-\frac{D_2}{2-D_2}\times\frac{K}{r}\times\frac{4(r-2\alpha-\sigma^2}{\gamma^2}\right)^{\frac{1}{2}}\\
&A_2^{\infty}=\frac{2}{2-D_2}\times\frac{K}{r}\times(x^*)^{-{D_2}}
\end{cases}
\end{equation}
\end{theorem}
\begin{prooftheorem}
The proof of the theorem follows immediately of the next two propositions. \qed
\end{prooftheorem}
\begin{prop}
\label{convergencia}
The solution of truncated problem (\ref{modified_problem}) converges to the solution of original one (\ref{optimal}).
\end{prop}
\begin{proofprop}
Notice that,  $X(s)\geq X^C(s)$, with probability 1 and thus the following inequalities hold: 
\begin{equation*}
E\left[\int_0^{\tau}e^{-r s}\Pi(X(s))ds\right]\geq E\left[\int_0^{\tau}e^{-r s}\Pi(X^C(s))ds\right]\quad\Rightarrow\quad V(x)\geq V_C(x).
\end{equation*}
Therefore, $V_C(x)\leq V(x)\quad, \forall C>0$ and consequently $\lim\sup_{C\to + \infty} V_C(x)\leq V(x)$.
In order to prove that $\lim\inf_{C \to + \infty} V_C(x)\geq V(x)$, we will show that $\lim\inf_{C\to + \infty}J_C(x,\tau)\geq J(x,\tau),\quad \forall \tau\in{\cal S}$. Fix $\tau\in{\cal S}$. Then,
\begin{align*}
\label{igualdade}
E&\left[\int_0^{\tau}e^{-r s}\Pi(X(s))ds\right]=E\left[\int_0^{\tau\wedge\nu^c}e^{-r s}\Pi(X^C(s))ds\right]+\\
+&E\left[\int_{\tau\wedge\nu^C}^{\tau}e^{-r s}\Pi(X(s))ds\right]=\\
=&E\left[\int_0^{\tau\wedge\nu^C}e^{-r s}\Pi(X^C(s))ds\right]+E\left[\int_{\tau\wedge\nu^C}^{\tau}e^{-r s}\Pi(0)ds\right]+\\
+&E\left[\int_{\tau\wedge\nu^C}^{\tau}e^{-r s}\big(\Pi(X(s))-\Pi(0)\big)ds\right]=\\
=&E\left[\int_0^{\tau}e^{-r s}\Pi(X^C(s))ds\right]+E\left[\int_{\tau\wedge\nu^C}^{\tau}e^{-r s}\big(\Pi(X(s))-\Pi(0)\big)ds\right],
\end{align*}
where $a\wedge b=\min(a,b)$. By construction, $X^C(s)\nearrow X(s)$ a.s, when $C\to +\infty$. Thus it follows trivially that $\Pi(X^C(s))\nearrow\Pi(X(s))$ a.s. Furthermore, $\tau\wedge\nu^C\nearrow\tau$, for the same reason. So we have that: 
$$
\int_{\tau\wedge\nu^C}^{\tau}e^{-r s}\big(\underbrace{\Pi(X(s))-\Pi(0)}_{\geq 0}\big)ds\searrow 0, \quad\text{a.s}
$$
and, as $\Pi(X(s))-\Pi(0)\geq 0$ a.s, it follows that:
$$
\int_{\tau\wedge\nu^C}^{\tau}e^{-r s}\big(\underbrace{\Pi(X(s))-\Pi(0)}_{\geq 0}\big)ds\leq\int_0^{+\infty}e^{-r s}\big(\Pi(X(s))-\Pi(0)\big)ds.
$$
Therefore, the dominated convergence theorem guarantees that:
$$
\lim_{C\to\infty} E\left[\int_{\tau\wedge\nu^C}^{\tau}e^{-r s}\big(\Pi(X(s))-\Pi(0)\big)ds\right]=0.
$$\qed
\end{proofprop}
Given that convergence is guaranteed, the solution of the original problem is $V(x)=\lim_{C\to +\infty}V_C(x)$. Indeed, we only need to obtain the solution of system given by (\ref{fitconditions}) when $C\to\infty$.
\begin{prop}
The solution of the system (\ref{fitconditions}) when $C\to +\infty$ is given by (\ref{solution_infty}).
\end{prop}
\begin{proofprop}
In order to prove the result, notice that (\ref{fitconditions}) is equivalent to
\begin{equation}
\label{modifiedsystem}
\begin{cases}
A_1C^{D_1}+A_2C^{D_2}+\frac{\gamma^2}{4(r-2\alpha-\sigma^2)}C^2-\frac{K}{r}&=0\\
A_1(x_C^*)^{D_1}+A_2(x_C^*)^{D_2}+\frac{\gamma^2}{4(r-2\alpha-\sigma^2)}(x_C^*)^2-\frac{K}{r}&=0\\
A_1D_1(x_C^*)^{D_1}+D_2A_2(x_C^*)^{D_2}+\frac{\gamma^2}{2(r-2\alpha-\sigma^2)}(x_C^*)^2&=0
\end{cases}
\end{equation}
and that dividing the first equation by $C^{D_1}$ we obtain
$A_1+\epsilon_1A_2+\epsilon_2=0$ where $\epsilon_1=C^{D_2-D_1}$, $\epsilon_2=BC^{2-D_1}-\frac{K}{r}C^{-D_1}$ and $B=\frac{\epsilon^2}{4(r-2\alpha-\sigma^2)}$. As $D_2<0$ and $D_1>2$, $\epsilon_1,\epsilon_2\to 0$ when $C\to +\infty$. Making $\epsilon_1,\epsilon_2=0$, it is possible to solve explicitly the system, and its solution is given by (\ref{solution_infty}).
In order to prove the result, we will use the Implicit Function Theorem.
Now let $\psi:\mathbb{R}^5\to\mathbb{R}^3$ be the $C^\infty$ function defined by:
\begin{align*}
&\psi(A_1,A_2,x,\epsilon_1,\epsilon_2)=\\
&\left(A_1+\epsilon_1 A_2+\epsilon_2, A_1x^{D_1}+A_2x^{D_2}+Bx^2-\frac{K}{r}, A_1D_1x^{D_1}+D_2A_2x^{D_2}+2Bx^2\right)
\end{align*}
The previous considerations guarantee that $\psi(A^{\infty}_1,A^{\infty}_2,x^*_{\infty},0,0)=0$, so, if we represent the matrix of the first derivatives of $\psi$ in the variables $A_1,A_2$ and $x$ by $D_{(A_1,A_2,x)}\psi$ , we will need to verify that $|D_{(A_1,A_2,x)}\psi(A^{\infty}_1,A^{\infty}_2,x^*_{\infty},0,0)|\neq 0$. As,
\begin{align*}
D_{(A_1,A_2,x)}&\psi(A^{\infty}_1,A^{\infty}_2,x^*_{\infty},0,0) =\\
  &\begin{pmatrix}
   1 & 0 & 0 \\
   (x^*_{\infty})^{D_1} & (x^*_{\infty})^{D_2} & D_2A_2(x^*_{\infty})^{D_2-1}+2B(x^*_{\infty}) \\
   D_1(x^*_{\infty})^{D_1} & D_2(x^*_{\infty})^{D_2} & D_2^2A_2(x^*_{\infty})^{D_2-1}+4B(x^*_{\infty})  \end{pmatrix},
\end{align*}
trivial calculations show that $|D_{(A_1,A_2,x)}\psi(A^{\infty}_1,A^{\infty}_2,x^*_{\infty},0,0)|=$\\$(4-2D_2)B(x^*_{\infty})^{D_2+1}\neq 0$. Then, the Implicit Function Theorem, guarantees that there are an open set $U\subset\mathbb{R}^3$ containing $(A^{\infty}_1,A^{\infty}_2,x^*_{\infty})$ and an open set $V\subset\mathbb{R}^2$ containing $(0,0)$ and a unique continuously differentiable function $\phi:V\to U$ such that: 
\begin{align*}
\{\psi&\left(\phi(\gamma_1,\gamma_2),\gamma_1,\gamma_2\right):(\gamma_1,\gamma_2)\in V\}=\\
&\{\psi\left(A^{\infty}_1,A^{\infty}_2,x^*_{\infty},\gamma_1,\gamma_2\right):(\gamma_1,\gamma_2)\in V,(A^{\infty}_1,A^{\infty}_2,x^*_{\infty})\in U\}.
\end{align*}
So, this guarantees that, when $(\gamma_1,\gamma_2)\to (0,0)$ the solution of \\$\psi(A_1,A_2,x,\gamma_1,\gamma_2)=0$ is given by (\ref{solution_infty}). \qed
\end{proofprop}

\section{Conclusion}

The way that we derive the solution for the value of a firm that faces the option to exit the market can also be used to solve other investment problems, as investment in new products, suspension of production, etc.

In fact in all the cases of investment decisions the problem turns out to be a free-boundary one. But to apply the verification theorem that we used in order to prove existence and uniqueness of classical solution we need a sufficient number of boundary conditions. When we consider the truncation of the original model to obtain the necessary boundary conditions. 

A possible  extension of the present  work is the introduction of suspension option in the production process. There we may assume that, because there are running costs, the firm may temporarly stop production, entering in a suspension mode. Production may be resumed latter, paying some cost of return. This new problem should be formalized as an impulse control problem, for which we pretend to obtain the associated HJB equation and its solution.

\section{Appendix}
\label{A}

Consider the unidimensional stochastic process, $\{X(t), t\geq 0\}$ defined on a complete filtered space $((\Omega, \{{\cal F}_t\}_{t\geq 0}, {\cal P})$ satisfying the following diffusion equation 
\begin{align}
&dX(t)=b(X(t))dt+\sigma(X(t))dW(t)\\
&X(0)=x
\end{align}
Moreover we assume that $b(.)$ and $\sigma(.)$ satisfy the usual conditions of measurability and adaptability, see for example \cite{oksendal}. Let $\tau$ be a Markov time with respect to the stochastic process $\{X(t), t\geq 0\}$, taking values in $[0,+\infty]$.
We define ${\cal S}$ as the state space of all Markov Times adapted to the filtration generated by the Brownian Motion defined above.

Suppose that the continuous function $g,h:\mathbb{R}\mapsto \mathbb{R}$ satisfy the following measurability condition:
\begin{equation}
\label{condition}
E\left[\int^{+\infty}_0e^{-r t}|g(X(t))|dt+e^{-r \tau}|h(X(\tau))|\right]<\infty.   
\end{equation} 
Let
\begin{equation}
\label{problem}
J(x,\tau)\equiv E\left[\int_0^{\tau}e^{-r s}g(X(s))ds+e^{-r \tau}h(X(\tau))\chi_{\tau<\infty}\right]
\end{equation}
then the optimal stopping time problem is to find a function $V$ such that  
$$V(x)=\sup_{\tau\in {\cal S}}   J(x,\tau).$$

In the following Theorem we provide a characterization of the solution of  (\ref{problem}) as a solution of a second order differential equation .  Its proof can be found, for instance, in \cite{oksendal}.
\begin{theorem}
\label{verification}
Let $\phi(.)$ be a function, with $\phi\in C^1(\mathbb{R})$, and with derivative absolutely continuous, such that $\phi$ satisfies the following condition
$$
\min\{r \phi(x)-b(x)\phi'(x)-\frac{\sigma^2(x)}{2}\phi''(x)- g(x), \phi(x)-h(x)\}=0 
$$
for almost every $x\in\mathbb{R}$. \end{theorem}

\bibliographystyle{elsarticle-harv}
\bibliography{myrefs}

\end{document}